\documentstyle[11pt]{amsart}

\def\demo{\noindent{\bf Proof. }}
\def\sqr#1#2{{\vcenter{\hrule height.#2pt
            \hbox{\vrule width.#2pt height#1pt \kern#1pt
                    \vrule width.#2pt}
            \hrule height.#2pt}}}
\def\square{\mathchoice\sqr64\sqr64\sqr{4}3\sqr{3}3}
\def\QED{\hfill$\square$}

\newtheorem{Theorem}{Theorem}[section]

\newtheorem{Corollary}[Theorem]{Corollary}
\newtheorem{Proposition}[Theorem]{Proposition}

\newtheorem{Remark}[Theorem]{Remark}
\newtheorem{Example}[Theorem]{Example}

\newtheorem{Question}[Theorem]{Question}

\begin{document}

\baselineskip=17pt

\title[Depth of associated graded rings via Hilbert
coefficients of ideals] {\bf Depth of associated graded rings
via \\ Hilbert coefficients of ideals}

\author[A. Corso, C. Polini and M.E. Rossi]
{Alberto Corso \and Claudia Polini \and Maria Evelina Rossi}

\address{Department of Mathematics, University of Kentucky,
Lexington, Kentucky 40506 - USA}
\email{corso@@ms.uky.edu}

\address{Department of Mathematics, University of Notre Dame,
Notre Dame, Indiana 46556 - USA} \email{cpolini@@nd.edu}

\address{Dipartimento di Matematica, Universit\`a di Genova, Via Dodecaneso 35,
16132 Genova - Italy} \email{rossim@@dima.unige.it}

\dedicatory{Dedicated to Wolmer V. Vasconcelos on the occasion of
his sixty fifth birthday}

\thanks{
AMS 2000 {\em Mathematics Subject Classification}. Primary 13A30,
13B21, 13D40; Secondary 13H10, 13H15. \newline\indent The first
author was partially supported by a Faculty Summer Research
Fellowship from the University of Kentucky. The second author
gratefully acknowledges partial support from the NSF and the AWM
Mentoring Travel Grant. The third author was partially supported
by the Italian MIUR}

\begin{abstract}
Given a local Cohen-Macaulay ring $(R, {\mathfrak m})$, we study
the interplay between the integral closedness -- or even the
normality -- of an ${\mathfrak m}$-primary $R$-ideal $I$ and
conditions on the Hilbert coefficients of $I$. We relate these
properties to the depth of the associated graded ring of $I$.
\end{abstract}

\vspace{-0.1in}

\maketitle

\vspace{-0.2in}

\section{Introduction}

\noindent Let $I$ be an ${\mathfrak m}$-primary ideal of a local
Cohen-Macaulay ring $(R, {\mathfrak m})$ of dimension $d>0$ and
with infinite residue field. The {\it Hilbert-Samuel function} of
$I$ is the numerical function that measures the growth of the
length of $R/I^n$ for all $n \geq 1$. For $n \gg 0$ this function
$\lambda(R/I^n)$ is a polynomial in $n$ of degree $d$
\[
\lambda(R/I^n) = e_0 {{n+d-1}\choose d} - e_1 {{n+d-2}\choose d-1}
+ \ldots + (-1)^d e_d,
\]
and  $e_0, e_1, \ldots, e_d$ are called the Hilbert coefficients of $I$.

It is well-known that $e_0 =\lambda(R/J)$ for any minimal
reduction $J$ of $I$ and that the integral closure $\overline{I}$
of $I$ can also be characterized as the largest ideal containing
$I$ with the same {\it multiplicity} $e_0$ \cite {NR}. More
generally, L.J. Ratliff and D.E. Rush introduced the ideal
$\widetilde{I}$, which turns out to be the largest ideal
containing $I$ with the same Hilbert coefficients as $I$
\cite{RR}. In particular one has the inclusions $I \subseteq
\widetilde I \subseteq \overline{I}$, where equalities hold if $I$
is integrally closed. A very useful technique -- that we also
exploit -- is to consider the generating functions of
$\lambda(R/\widetilde {I^n})$ or $\lambda(R/\overline {I^n})$
instead of the one of $\lambda(R/{I^n})$: They clearly coincide if
$I$ is normal $($that is, all powers of $I$ are integrally
closed$)$.

Little is known about the higher Hilbert  coefficients of $I$,
unless we are in presence of good depth properties of the {\it
associated graded ring} ${\rm gr}_I(R)=
\displaystyle\bigoplus_{n\geq 0}I^n/I^{n+1} $ of $ I. $  For
example, if the depth of ${\rm gr}_I(R)$ is at least $d-1$ then
all the Hilbert  coefficients of $I$ are positive \cite{M}.
Conversely, numerical information on the $e_i$'s has been used to
obtain information on the depth of ${\rm gr}_I(R)$. For instance,
D.G. Northcott showed that the bound $\lambda(R/I) \geq e_0 - e_1$
always holds \cite{No}. Later C. Huneke \cite{H} and A. Ooishi
\cite{O} showed that the equality $\lambda(R/I) = e_0 - e_1$ holds
if and only if $I^2=JI$.  In particular, ${\rm gr}_I(R)$ is
Cohen-Macaulay.

Translating information from the Hilbert  coefficients of
$I$ into good depth properties of ${\rm gr}_I(R)$ has also been a
constant theme in the work of J.D. Sally
\cite{S,S1,S2,S3,S4,S6,S5}. The most recent results along this
line of investigation can be found in
\cite{CPVP,E,ERV,EV,Huc,HM,Itoh,P,Ro,RV,V,Vaz,W}. The
general philosophy is that an `extremal' behavior of some of the
$e_i$'s controls the depth of the associated graded ring of $I$,
or of some of its powers, and at the same time forces its
Hilbert-Samuel function. We remark that these results are somewhat
unexpected since the Hilbert coefficients give
asymptotic information on the Hilbert-Samuel function.

It is clear that $e_0$ and $e_1$ are positive integers. As far as
the higher Hilbert coefficients of $I$ are concerned, it is
a famous result of M. Narita that $e_2 \geq 0$ \cite{Nar}. In this
case the minimal value for $e_2$ does not imply the
Cohen-Macaulayness of ${\rm gr}_{I }(R)$. In the very same paper,
he also showed that if $d=2,$ then $e_2=0$ if and only if $I^n$ has
reduction number one for some $n \gg 0$. In particular, ${\rm
gr}_{I^n}(R)$ is Cohen-Macaulay. Examples show that the result
cannot be extended to higher dimension. In \cite{CPVP} an
elementary proof of the positivity of $e_2 $ has been given by
using the structure of the so-called Sally module $S_J(I)$.

Unfortunately, the well-behavior of the Hilbert
coefficients stops with $e_2$. Indeed, in \cite{Nar} M. Narita
showed that it is possible for $e_3$ to be negative. However, a
remarkable result of S. Itoh says that if $I$ is a normal ideal
then $e_3 \geq 0$ \cite{I}. A recent proof of this result was
given by S. Huckaba and C. Huneke in \cite{HH}. In general, it
seems that the integral closedness $($or the normality$)$ of the
ideal $I$ yields non trivial consequences on the Hilbert
coefficients of $I$ and, ultimately, on ${\rm depth}\,{\rm
gr}_I(R)$.

To be more specific, our goal is to characterize a sufficiently
high depth of the associated graded ring of $I$ in terms of
conditions on the first Hilbert coefficients, and in
particular on $e_2$ and $e_3$. Our approach is to study the
interplay between the integral closedness $($or the normality$)$
of the ideal $I$ and $($upper or lower$)$ bounds on the
Hilbert  coefficients of $I$ and relate it to the depth of
the corresponding associated graded ring of $I$. Among our tools,
we make systematic use of the standard technique of modding out a
superficial sequence in order to decrease the dimension of the
ring. This explains why some of our results are formulated for
rings of small dimension.

We first establish in Theorem~\ref{thm0} a general upper bound on
the second Hilbert  coefficient which is reminiscent of a similar
bound on the first Hilbert coefficient due to S. Huckaba and T.
Marley \cite{HM} and M. Vaz Pinto \cite{Vaz}. Namely, we show that
$e_2 \leq \displaystyle \sum_{n \geq 1} n \lambda(I^{n+1}/JI^n)$,
for any minimal reduction $J$ of $I$. Furthermore, the upper bound
is attained if and only if ${\rm depth}\, {\rm gr}_I(R) \geq d-1$.
Next, we characterize in Theorem~\ref{thm0.5} the depth of the
associated graded ring of ideals whose second Hilbert coefficient
has value `close enough' to the upper bound established in
Theorem~\ref{thm0}. That is, the condition $e_2 \geq \displaystyle
\sum_{n \geq 1} n \lambda(I^{n+1}/JI^n) - 2$ implies that ${\rm
depth}\, {\rm gr}_I(R) \geq d-2$. Noteworthy is the fact that the
same conclusion holds whenever $I$ is an integrally closed ideal
satisfying the less restrictive bound $e_2 \geq \displaystyle
\sum_{n \geq 1} n \lambda(I^{n+1}/JI^n)-4$. Still in the case of
integrally closed ideals we improve Narita's positivity result on
$e_2$. Indeed, in Theorem~\ref{thm1} we show that for any
integrally closed ideal $I$ one has   $e_2 \geq \lambda(I^2/JI)$
for any minimal reduction $J$ of $I$. The interesting fact is that
equality holds in the previous formula if and only if $I^3=JI^2$
if and only if $\lambda(R/I) = e_0-e_1+\lambda(I^2/JI)$. In this
case ${\rm gr}_I(R)$ is Cohen-Macaulay and the Hilbert-Samuel
function is forced. This result fully generalizes the ones of Itoh
\cite{Itoh}, who handled the cases $e_2 \leq 2$.

In \cite{S5} J.D. Sally proved that if $d=2$ then $e_2 \ge e_1 -
e_0 +\lambda(R/\widetilde I)$. Starting from Sally's inequality,
S. Itoh proved that if $d \ge 1 $ and $I $ is integrally closed
then $e_2 \ge e_1 - e_0 +\lambda(R/I)$ \cite{Itoh}. It was not
known which conditions are forced on the ideal $I$ if the equality
$e_2 = e_1 - e_0 +\lambda(R/I)$ holds whenever $I$ is integrally
closed. In the particular case of the maximal ideal ${\mathfrak
m}$ of $R$, G. Valla had conjectured in \cite{Val} that the
reduction number of ${\mathfrak m}$ is at least two, which in turn
would imply the Cohen-Macaulayness of ${\rm gr}_{{\mathfrak
m}}(R)$. A recent example of H.-S. Wang shows, though, that
Valla's conjecture is false in general. In Theorem~\ref{valla} we
prove Valla's conjecture for normal ideals. In particular, if $I$
is normal and \ $\lambda(R/I)=e_0-e_1+e_2$ \ then $e_2=
\lambda(I^2/JI)$ and $I^3=JI^2$ for some minimal reduction $J$ of
$I$. In particular, ${\rm gr}_I(R)$ is Cohen-Macaulay and the
Hilbert function is known. The key to the result is a theorem of
Itoh on the normalized Hilbert  coefficients of ideals generated
by a system of parameters.

As far as the third Hilbert  coefficient of $I$ is
concerned, our first result in Section~4 is a generalization of
Itoh's result on the positivity of $e_3$ in case $d=3.$  The
thrust of our calculation is to replace the normality assumption
on $I$ with the weaker requirement of the integral closedness of
$I^n$ for some large $n$ $($see Theorem~\ref{thm3}$)$. The proof
reduces to comparing the Hilbert  coefficients of $I$ and
those of a large power of $I.$
Combining this result with Theorem~\ref{valla} we are able to
characterize when $e_3=0$ for asymptotically normal ideals. If
this is the case, then for $n \gg 0$ we have that $I^n$ has
reduction number at most two, which in turn yields that ${\rm
gr}_{I^n}(R)$ is Cohen-Macaulay. This result is reminiscent of
Narita's characterization of $e_2=0$ when $d=2$.

\section{Preliminaries}

\noindent Thus far we have described the Hilbert-Samuel function
  associated
with the $I$-adic filtration ${\mathcal F} =\{ I^n \}_{n\geq 0}$.
It is important to observe that the theory also applies to other
filtrations of ideals of $R$: The so-called Hilbert filtrations
(see \cite{HM} and \cite{GR}). Let $(R, {\mathfrak m})$ be a local
ring of dimension
$d$. A filtration of $R$-ideals ${\mathcal F}= \{ F_n \}_{n\geq
0}$ is called an {\it Hilbert filtration} if $F_1$ is an
${\mathfrak m}$-primary ideal and the Rees algebra ${\mathcal
R}({\mathcal F})$ is a finite ${\mathcal R}(F_1)$-module. As in
the case of the $I$-adic filtration where $I$ is an ${\mathfrak
m}$-primary ideal, we can define the {\it Hilbert-Samuel function}
of ${\mathcal F}$ to be $\lambda(R/F_n)$. For $n\gg 0$ one also
has that $\lambda(R/F_n)$ is a polynomial in $n$ of degree $d$
\[
\lambda(R/F_n) = \sum_{j=0}^d (-1)^j e_j({\mathcal F})
{n+d-j-1\choose d-j},
\]
where the $e_j ({\mathcal F})$'s are called the {\it
Hilbert  coefficients} of ${\mathcal F}$.

Another related object is the {\it Hilbert series} of ${\mathcal
F}$, which is defined as
\[
P_{\mathcal F}(t) = \sum_{n\geq 1} \lambda(F_{n-1}/F_n) t^{n-1}.
\]
The numerical function $\lambda(F_{n-1}/F_n)$ is called {\it
Hilbert function} with respect to the filtration ${\mathcal F}$.
It is well known that there exists a unique polynomial
$f_{\mathcal F}(t) \in {\mathbb Z}[t]$, called the {\it
$h$-polynomial} of ${\mathcal F}$, with degree $s({\mathcal F})$,
$f_{\mathcal F}(1) \not=0$ and such that
\[
P_{\mathcal F}(t) = \frac{f_{\mathcal F}(t)}{(1-t)^d}.
\]
We recall that $e_j({\mathcal F}) = f^{(j)}_{\mathcal F}(1)/j!$,
where $f^{(j)}_{\mathcal F}(t)$ denotes the $j$-th formal
derivative of $f_{\mathcal F}(t)$, and we also point out that it
is useful to consider the Hilbert coefficients $e_j
({\mathcal F})$ even when $j > d$. Finally, we denote by ${\rm
gr}_{{\mathcal F}}(R)= \displaystyle\bigoplus_{n\ge 0}F_n/F_{n+1}$
the associated graded ring with respect to ${\mathcal F}$.

If $I$ is an ${\mathfrak m}$-primary ideal and ${\mathcal F} = \{
I^n \}_{n\geq 0}$ is the usual $I$-adic filtration we write
$e_j(I)$ instead of $e_j ({\mathcal F})$ or simply $e_j$ if there
is not confusion on the ideal under consideration and we denote by
${\rm gr}_I(R)$ the corresponding associated graded ring.  Of
particular interest is the filtration ${\mathcal F} = \{
\overline{I^n} \}_{n\geq 0}$ given by the integral closure of the
powers of an ${\mathfrak m}$-primary ideal $I$ of an analytically
unramified local ring. It is customary to denote with
$\overline{e}_0(I), \overline{e}_1(I), \ldots, \overline{e}_d(I) $
the Hilbert coefficients with respect to this filtration
${\mathcal F}$. We also recall that if $J$ is a  reduction of $I$
$($that is $J \subset I$ and $I^{n+1} = JI^n$ for some integer
$n)$, then $\overline{J^n}= \overline{I^n}$ for every $n$ because
$J^n$ is a reduction of $I^n. $ It follows that if $I$ is normal
then $e_j(I) = \overline{e}_j(J)$. Another crucial example is the
Ratliff-Rush filtration of the powers of an ${\mathfrak
m}$-primary ideal $I$ of a local ring. We recall that \cite{RR}
\[
\widetilde{I^n} = \bigcup_{k \ge 1}I^{k+n}: I^k.
\]
If $I$ contains a non zero divisor then $\widetilde {I^n}= I^n$
for $n \gg 0$ and hence ${\mathcal F} = \{ \widetilde {I^n}
\}_{n\geq 0}$ is an Hilbert filtration. In particular $e_j
({\mathcal F})=e_j(I) $ for $j=0, \dots,d.$

The advantage of considering ${\rm gr}_{{\mathcal F}}(R)$ with the
above filtrations rather than ${\rm gr}_I(R) $ is that they are
graded rings with positive depth. Unfortunately, they are not
standard algebras.

A classical technique for studying the Hilbert  coefficients of
any filtration ${\mathcal F}$ is to reduce the dimension of the
ring by modding out a superficial sequence for ${\mathcal F}$. We
recall that an element $x \in F_1$ is called a {\it superficial
element} for ${\mathcal F}$ if there exists an integer $c$ such
that $(F_n \colon x) \cap F_c =F_{n-1}$ for all $n>c$. A sequence
$x_1, \ldots, x_k$ is then called a {\it superficial sequence} for
${\mathcal F}$ if $x_1$ is superficial for ${\mathcal F}$ and
$x_i$ is superficial for the quotient filtration ${\mathcal
F}/(x_1, \ldots, x_{i-1})= \{ F_n + (x_1, \ldots, x_{i-1})/(x_1,
\ldots, x_{i-1})\}$ for $2 \leq i \leq k$. Now, if ${\rm grade}\,
F_1 \geq k$ and $x_1, \ldots, x_k$ is a superficial sequence for
${\mathcal F}$ it can be showed that $e_j({\mathcal
F})=e_j({\mathfrak F})$, for $0\leq j \leq d-k$, where ${\mathfrak
F}={\mathcal F}/(x_1, \ldots, x_k)= \{ F_n + (x_1, \ldots,
x_k)/(x_1, \ldots, x_k)\}$ (see for instance \cite{HM}).

Sometimes, though, it may not be possible to lift certain results
back to the original ideal.  That's the reason why Narita's
characterization of $e_2=0$, that we mentioned in the
introduction, only works in dimension two.

In \cite{I}, Itoh systematically used the reduction method in the
case of normal ideals. He proved that if $I$ is integrally closed,
then there exists a superficial element $x $ for $I$ such that
$I/(x) $ is still integrally closed. He had to be particularly
careful as the normality is not preserved by going modulo a
superficial sequence. If $I$ is normal and $x$ is a superficial
element for $I$, Morales \cite{Mo} and independently Itoh
\cite{I} proved that $I/(x)$ is still normal if and only
if depth ${\rm gr}_I(R) \ge 2$. With some clever techniques Itoh
could show, though, that all sufficiently large powers of $I/(x)$
remain integrally closed. He then managed to prove his result on
the positivity of $e_3$ using a great deal of local cohomology
calculations.

\section{Results on the second Hilbert  coefficient}

\noindent Our first result is an upper bound on $e_2$ which is
reminiscent of the one on $e_1$ established by Huckaba and Marley
\cite[4.7]{HM} and by Vaz Pinto \cite[1.1]{Vaz}, which characterizes
when the depth of the associated graded ring is at least $d-1$.

\begin{Theorem}\label{thm0}
Let $(R, {\mathfrak m})$ be a local Cohen-Macaulay ring of
dimension $d \geq 1$ and let $I$ be an ${\mathfrak m}$-primary
ideal of $R$. Then
\[
e_2 \le \displaystyle\sum_{n\ge 1} n\lambda(I^{n+1}/JI^n),
\]
for any minimal reduction $J$ of $I$. Furthermore, equality holds
for some minimal reduction $J$ of $I$ if and only if \/ ${\rm
depth}\, {\rm gr}_I(R) \geq d-1$.
\end{Theorem}
\demo If $d=1$ the result follows from \cite[1.9]{GR}. Let us
assume then $d \geq 2$. Let $J$ be a minimal reduction of $I $ and
let $x_1, \dots, x_{d-1}$ be a superficial sequence for $I$
contained in $J$. Let $H$ and $K$ denote the ideals $I/(x_1,
\dots,x_{d-2})$ and $I/(x_1, \dots,x_{d-1})$, respectively. By
\cite[1.2$({\it a})$,$({\it b})$]{ERV} and \cite[1.9]{GR} we get
$e_2(I)=e_2(H) \le e_2(K) = \displaystyle\sum_{n\ge 1}
n\lambda(K^{n+1}/JK^n) \le \displaystyle\sum_{n\ge 1}
n\lambda(I^{n+1}/JI^n)$, which establishes the desired inequality.

If depth ${\rm gr}_I(R) \geq d-1$ the equality follows from
\cite[1.9]{GR}. Conversely, if  equality holds one has that
$e_2(H) = e_2(K)$, which in turn forces ${\rm depth}\, {\rm
gr}_H(R/(x_1, \dots,x_{d-2})) \geq 1$ by \cite[1.2$({\it
c})$]{ERV}. Hence by \cite[2.2]{HM} we conclude that ${\rm
depth}\, {\rm gr}_I(R) \geq d-1$. \QED

\medskip

The following example, due to Huckaba and Huneke \cite[3.12]{HH},
provides an instance in which the bound in Theorem~\ref{thm0} is
attained. This example will also play a role in the next section.

\begin{Example}\label{ex-HuHu}{\rm
Let $k$ be a field of characteristic $\not= 3$ and set $R =
k[\![X,Y,Z]\!]$, where $X,Y,Z$ are indeterminates. Let $N = (X^4,
X(Y^3+Z^3), Y(Y^3+Z^3), Z(Y^3+Z^3))$ and set $I = N + {\mathfrak
m}^5$, where ${\mathfrak m}$ is the maximal ideal of $R$. The
ideal $I$ is a normal ${\mathfrak m}$-primary ideal whose
associated graded  ring ${\rm gr}_I(R)$ has depth $d-1$, where
$d(=3)$ is the dimension of $R$. We checked that
\[
P_I(t)=\frac{31 + 43t +t^2 + t^3}{(1-t)^3},
\]
thus yielding  $ e_2=4$. Moreover, we also checked that
$\lambda(I^2/JI)=2$, $\lambda(I^3/JI^2)=1$ and $I^4=JI^3$, for any
minimal reduction $J$ of $I$. Hence the bound in
Theorem~\ref{thm0} is sharp.

Since $I$ is, in particular, integrally closed we have that $J
\cap I^2 = JI$ by \cite[2.1]{H} and \cite[1]{It}. Thus ${\rm
depth}\, {\rm gr}_I(R) \geq 2$ also follows from the main result of
\cite{CPVP,E,Huc2,R}. }
\end{Example}

\medskip

Theorem~\ref{thm0.5} below deals with the depth property of the
associated graded ring of ideals whose second Hilbert
coefficient is close enough to the upper bound established in
Theorem~\ref{thm0}.

\begin{Theorem}\label{thm0.5}
Let $(R, {\mathfrak  m})$ be a local Cohen-Macaulay ring of
dimension $d \geq 1$. Let $I$ be an ${\mathfrak m}$-primary ideal
and let $J$ denote a minimal reduction of $I$. If any of the
following conditions holds
\begin{itemize}
\item[$({\it a})$]
$e_2 \geq \displaystyle\sum_{n \geq 1} n \lambda(I^{n+1}/JI^n)-2;$

\item[$({\it b})$]
$I$ is an integrally closed ideal and $e_2 \geq
\displaystyle\sum_{n \geq 1} n \lambda(I^{n+1}/JI^n)-4,$
\end{itemize}
then ${\rm depth}\, {\rm gr}_I(R) \geq d-2$.
\end{Theorem}
\demo Throughout the proof we use the same notation as in the
proof of Theorem~\ref{thm0}.  By Theorem~\ref{thm0} we may assume
that ${\rm depth}\, {\rm gr}_I(R) < d-1$, which implies that
$e_2(H) < e_2(K)$ and
\[
\sum_{n\geq 1} n \lambda(K^{n+1}/JK^n) < \sum_{n \geq 1} n
\lambda(H^{n+1}/JH^n) \le \sum_{n \geq 1} n
\lambda(I^{n+1}/JI^n).
\]
Indeed $e_2(H) =e_2(K)$ implies ${\rm depth}\, {\rm gr}_H(R/(x_1,
\dots, x_{d-2})) >0$ by \cite[1.2$({\it c})$]{ERV} and hence ${\rm
depth} \, {\rm gr}_I(R) \geq d-1$  by \cite[2.2]{HM}.
If
\[
\sum_{n\geq 1} n \lambda(K^{n+1}/JK^n) = \sum_{n \geq 1} n
\lambda(H^{n+1}/JH^n),
\]
then $\lambda(K^{n+1}/JK^n) = \lambda(H^{n+1}/JH^n)$ for every $n$
so that
\[
e_1(H)=e_1(K)=\sum_{n\geq 0}  \lambda(K^{n+1}/JK^n) = \sum_{n \geq 0}
\lambda(H^{n+1}/JH^n),
\]
hence  again ${\rm depth} \, {\rm gr}_H(R/(x_1, \dots, x_{d-2}))
>0$ by \cite[4.7$({\it b})$]{HM} and, as before, $ {\rm depth}\  {\rm
gr}_I(R) \geq d-1 $.

Let us assume that $({\it a})$ holds. We have
\begin{eqnarray*}
\sum_{n \geq 1} n \lambda(I^{n+1}/JI^n) - 2 & \leq & e_2(I) =
e_2(H) \\
& \leq & e_2(K)-1 \\
& = & \sum_{n \geq 1} n \lambda(K^{n+1}/JK^n)-1 \\
& \leq & \sum_{n\geq 1} n \lambda(I^{n+1}/JI^n)-2.
\end{eqnarray*}
Hence we obtain that
\[
\sum_{n\geq 1} n \lambda(K^{n+1}/JK^n) = \sum_{n \geq 1} n
\lambda(I^{n+1}/JI^n)-1,
\]
which implies that $\lambda(K^2/JK)=\lambda(I^2/JI)-1$ and
$\lambda(K^{n+1}/JK^n)=\lambda(I^{n+1}/JI^n)$ for all $n \geq 2$.
Hence
\[
e_1(I) = e_1(K) = \sum_{n\geq 0} \lambda(K^{n+1}/JK^n) = \sum_{n
\geq 0} \lambda(I^{n+1}/JI^n)-1,
\]
from which it follows that ${\rm depth}\, {\rm gr}_I(R) \geq d -
2$ by \cite {P} and \cite {W}.

Let us assume now that $({\it b})$ holds. Since $I $ is integrally
closed $ I^2 \cap J= JI $ by \cite[4.7$({\it b} )$]{H} and
\cite[1]{It}, hence
$\lambda(H^2/JH)=\lambda(K^2/JK)=\lambda(I^2/JI)$. We claim that
our assumption on $e_2$ forces
$\lambda(H^{n+1}/JH^n)=\lambda(I^{n+1}/JI^n)$ for all $n \geq 0$
and $\displaystyle\sum_{n \geq 0} \lambda(K^{n+1}/JK^n) =
\displaystyle\sum_{n\geq 0} \lambda(I^{n+1}/JI^n) - 1$. Thus we
conclude, as before, that
\[
e_1(I) = \sum_{n \geq 0} \lambda(I^{n+1}/JI^n)-1,
\]
which forces ${\rm depth}\, {\rm gr}_I(R) \geq d-2$.

Notice that if $\lambda(H^{n+1}/JH^n) \not= \lambda(I^{n+1}/JI^n)$
for some $n\  (\geq 2) $ then
\begin{eqnarray*}
\sum_{n \geq 1} n \lambda(H^{n+1}/JH^n)-2 & \leq & \sum_{n \geq 1}
n \lambda(I^{n+1}/JI^n) -4 \\
& \leq & e_2(I) = e_2(H) \\
& \leq & e_2(K)-1 \\
& = & \sum_{n \geq 1} n\lambda(K^{n+1}/JK^n) -1 \\
& = & \lambda(H^2/JH) + \sum_{n \geq 2} n \lambda(K^{n+1}/JK^n)-1
\\
& \leq & \lambda(H^2/JH) + \sum_{n\geq 2} n \lambda(H^{n+1}/JH^n)
-3
\end{eqnarray*}
which is impossible. Hence we may assume that
$\lambda(I^{n+1}/JI^n)=\lambda(H^{n+1}/JH^n)$ for all $n. $ Since $
e_2(H) \le e_2(K)-1, $  we have to
consider two cases.

If $e_2(H)=e_2(K)-1$, by \cite[1.2$({\it b})$]{ERV} we have that
$\displaystyle\sum_{n  \geq 1} \lambda((H^{n+1} \colon
x_{d-1})/H^n) =1$. This implies $($by using, for example, the
exact sequence in the proof of 1.7 \cite{Ro}) that
$\displaystyle\sum_{n \geq 1} \lambda(H^{n+1}/JH^n) =
\displaystyle\sum_{n  \geq 1} \lambda(K^{n+1}/JK^n)+1$. This
proves our claim.

If $e_2(H) \leq e_2(K) - 2$,  then we have
\begin{eqnarray*}
   e_2(H)  & \leq  & e_2(K) - 2 \\
& = & \sum_{n\geq 1} n \lambda(K^{n+1}/JK^n) - 2 \\
& \leq & \sum_{n \geq 1} n \lambda(H^{n+1}/JH^n) - 4,
\end{eqnarray*}
which implies $\displaystyle\sum_{n\geq 2} n \lambda(K^{n+1}/JK^n)
= \displaystyle\sum_{n \geq 2} n \lambda(H^{n+1}/JH^n) - 2$. This
proves again our claim. \QED

\medskip

\begin{Remark}{\rm
From the proof of Theorem~\ref{thm0.5} we conclude that $e_2(I)$
can never be equal to $\displaystyle\sum_{n\geq 1} n
\lambda(I^{n+1}/JI^n)-1$. Moreover, if $I$ is an integrally closed
ideal we also conclude that it can be neither $\displaystyle
\sum_{n\geq 1} n \lambda(I^{n+1}/JI^n)-1$ nor
$\displaystyle\sum_{n \geq 1} n \lambda(I^{n+1}/JI^n)-2$.}
\end{Remark}

\medskip

We illustrate Theorem~\ref{thm0.5} with the following example
which has been slightly modified from one suggested to us by Wang.

\begin{Example}\label{Wang} {\rm
Let $R$ be the three dimensional local Cohen-Macaulay ring
\[
k[\![X,Y,Z,U,V,W]\!]/(Z^2,ZU,ZV,UV,YZ-U^3,XZ-V^3),
\]
with $k$ a field and $X,Y,Z,U,V,W$ indeterminates. Let
$x,y,z,u,v,w$ denote the corresponding images of $X,Y,Z,U,V,W$ in
$R$. One has that the associated graded ring ${\rm gr}_{{\mathfrak
m}}(R)$ of ${\mathfrak m}=(x,y,z,u,v,w)$ has depth $d-2$, where
$d(=3)$ is the dimension of $R$. Indeed, we checked that
\[
P_{{\mathfrak m}}(t)= \frac{1 + 3t + 3t^3 - t^4}{(1-t)^3},
\]
so that $e_2=3$. Moreover $\lambda({\mathfrak m}^2/J{\mathfrak
m})=2$, $\lambda({\mathfrak m}^3/J{\mathfrak m}^2)=2$ and
${\mathfrak m}^4=J{\mathfrak m}^3$, where $J=(x,y,w)$. Thus
Theorem~\ref{thm0.5} applies. }
\end{Example}

Next, we present an improvement of Narita's positivity result on
$e_2$, which holds for any integrally closed ideal. We give a more
concrete lower bound and we characterize the integrally closed
ideals for which the minimal value of $e_2$ is attained.

\begin{Theorem}\label{thm1}
Let $(R, {\mathfrak m})$ be a local Cohen-Macaulay ring of
dimension $d \geq 1$. Let $I$ be an ${\mathfrak m}$-primary
integrally closed ideal. Then
\[
e_2 \geq \lambda(I^2/JI),
\]
where $J$ is any minimal reduction of $I$. In addition, the
following conditions are equivalent:
\begin{itemize}
\item[$(${\it a}$)$]
$e_2 = \lambda(I^2/JI)$;

\item[$(${\it b}$)$]
$I^3=JI^2$;

\item[$(${\it c}$)$]
$\lambda(R/I) = e_0-e_1+\lambda(I^2/JI)$.
\end{itemize}
Moreover, if any of the previous equivalent conditions holds then
${\rm gr}_I(R)$ is Cohen-Macaulay and
\[
P_{I}(t) = \frac{\lambda(R/I) + ( e_0 - \lambda(R/I) -
\lambda(I^2/JI))t + \lambda(I^2/JI) t^2 }{(1-t)^d}.
\]
\end{Theorem}
\demo Let $J$ be a minimal reduction of $I$. By \cite[4.7$({\it
a})$]{HM} we have the inequality
\[
e_1 \geq \displaystyle\sum_{n\ge 0}  \lambda(I^{n+1}/J\cap I^n)=
e_0 - \lambda(R/I) + \displaystyle\sum_{n\ge 1}
\lambda(I^{n+1}/J\cap I^n).
\]
Since $I$ is integrally closed, by \cite[12]{Itoh} we also have
that $e_2 \geq e_1 - e_0 + \lambda(R/I)$. If we now take into
account the above inequality on $e_1$ we conclude that
\[
e_2 \geq \displaystyle\sum_{n\ge 1}  \lambda(I^{n+1}/J\cap I^n).
\]
On the other hand, $I$ being integrally closed implies that $J
\cap I^2 = JI$ by \cite[4.7$({\it b})$]{H} and \cite[1]{It}. Hence
we have that
\[
e_2 \geq \lambda(I^{2}/J I) +\displaystyle\sum_{n\ge 2}
\lambda(I^{n+1}/J\cap I^{n+1})\geq \lambda(I^2/JI),
\]
which is the desired inequality.

Let us prove the equivalences. If $e_2 = \lambda(I^{2}/J I)$, then
for every $n \ge 2$ we have that $\lambda(I^{n+1}/J\cap
I^{n+1})=0$ and $e_1= \displaystyle\sum_{n\ge 0}
\lambda(I^{n+1}/J\cap I^n)= e_0 - \lambda(R/I) + \lambda(I^2/JI)$.
This proves $({\it c})$. If $({\it c})$ holds, by \cite[4.7$({\it
a})$]{HM} we have that ${\rm gr}_I(R)$ is Cohen-Macaulay. In
particular, we obtain that $J \cap I^{n+1} = JI^n$ for every $n$.
Hence $I^{n+1}=J\cap I^{n+1}= JI^n$ for $n \geq 2$. This yields
$({\it b})$. Suppose now that $({\it b})$ holds. Then we have that
$\displaystyle\sum_{n\ge 2} n\lambda(I^{n+1}/JI^n)=0 $. Now
Theorem \ref{thm0} also gives us the upper bound $e_2 \leq
\lambda(I^2/JI)$, so that $({\it a})$ follows.

As far as the Hilbert series is concerned, since ${\rm gr}_I(R)$
is Cohen-Macaulay it follows that $P_I(t)= \displaystyle
\frac{P_{I/J}(t)}{(1-t)^d}$. In particular, $P_{I/J}(t)$ is a
polynomial of degree $2$ because $I^3 \subseteq J$. If we write
$P_{I/J}(t)= h_0+h_1t +h_2t^2$, then we necessarily conclude that
$h_0=\lambda(R/I)$ and $h_2=e_2= \lambda(I^2/JI)$. \QED

\medskip

\begin{Remark}{\rm
We observe that the lower bound on $e_2$ given in
Theorem~\ref{thm1} is well defined, as $\lambda(I^2/JI)$ is always
independent of the minimal reduction $J$ of $I$ \cite{Va}. Also,
Theorem~\ref{thm1} recovers previous results by Itoh, who treated
the cases in which $e_2=0,1,2$: In such instances ${\rm gr}_I(R)$
always turns out to be Cohen-Macaulay \cite[5,6,7]{Itoh}. In
addition to fully treating the general case, we also describe the
Hilbert series of $I$.

We point out that if $e_2=3$ then ${\rm gr}_I(R) $ is not
necessarily Cohen-Macaulay even if $I$ is the maximal ideal of a
local Cohen-Macaulay ring (see Example~\ref{Wang2}).

Finally, we observe that in Theorem~\ref{thm1} the assumption  on
the ideal $I$ being `integrally closed' cannot be weakened. The
following example shows that $e_2=0$ does not imply the
Cohen-Macaulayness of ${\rm gr}_I(R)$. }
\end{Remark}

\begin{Example}\label{e2=0}{\rm
Let $R$ be the three dimensional regular local ring
$k[\![X,Y,Z]\!]$, with $k$ a field and $X,Y,Z$ indeterminates. The
ideal $I = (X^2-Y^2, Y^2-Z^2, XY, XZ, YZ)$ is not integrally
closed and
\[
P_I(t)=\frac{5 + 6t^2 - 4t^3 + t^4}{(1-t)^3}.
\]
In particular, $e_2=0 $ and ${\rm gr}_I(R)$ has depth zero. In
fact, computing $P_{I/(XY)}(t) $ we can see that  $XY$ is a
superficial element for $I $ whose initial form is a zero-divisor
in ${\rm gr}_I(R)$.}
\end{Example}

By using the techniques of this paper, we can also give here a
short proof of a result of Narita who characterized $e_2=0$ for
any ${\mathfrak m}$-primary ideal of a two dimensional local
Cohen-Macaulay ring.

\begin{Proposition}\label{Narita}
Let $(R, {\mathfrak m})$ be a local Cohen-Macaulay ring of
dimension two and let $I$ be an ${\mathfrak m}$-primary ideal.
Then $e_2=0$ if and only if $I^n$ has reduction number one for
some positive integer $n$. Under these circumstances then ${\rm
gr}_{I^n}(R)$ is Cohen-Macaulay.
\end{Proposition}
\demo We first recall that $e_2=e_2(I^m)$ for every positive
integer $m$. Assume $e_2=0$ and let $n$ be an integer such that
$\widetilde{I^n}=I^n$. By \cite[2.5]{S5}, $0= e_2(I^n) \ge
e_1(I^n) -e_0(I^n) + \lambda(R/\widetilde{I^n}) = e_1(I^n)
-e_0(I^n) + \lambda(R/ {I^n})$. Hence $e_1(I^n) - e_0(I^n) +
\lambda(R/ {I^n})=0$ because it cannot be negative by
Northcott's inequality. The result follows now by C. Huneke
\cite[2.1]{H} and A. Ooishi \cite[3.3]{O}. For the converse, if
$I^n$ has reduction number one for some $n$, then $e_2(I^n) = 0$
for example by \cite [2.4]{GR}. In particular $e_2(I)=e_2(I^n)=0$.
It is clear that if $I^n$ has reduction number one, then ${\rm
gr}_{I^n}(R)$ is Cohen-Macaulay (see \cite{VV}). \QED

\medskip

We remark that Narita's result cannot be extended to a local
Cohen-Macaulay ring of dimension $>2$. The ideal $I$ described in
Example~\ref{e2=0} satisfies $e_2=0$, however $I^m$ has not
reduction number one for every $m$. In fact, it is enough to
remark that $I$ has not reduction number one $({\rm gr}_{I}(R)$ is
not Cohen-Macaulay$)$ and $I^m=(X,Y,Z)^{2m}$ for $m >1$ which has
reduction number two.

\medskip

In \cite[12]{Itoh} Itoh  showed that if $I$ is an integrally
closed ideal then \ $e_2 \geq e_1 - e_0 + \lambda(R/I)$. Later, it
has been conjectured by Valla \cite[6.20]{Val} that if the
equality $e_2= e_1 - e_0 + \lambda(R/I)$ holds in the case in
which $I$ is the maximal ideal ${\mathfrak m}$ of $R$ then the
associated graded ring ${\rm gr}_{\mathfrak m}(R)$ is
Cohen-Macaulay. Unfortunately, the following example given by Wang
shows that the conjecture is, in general, false.

\begin {Example} \label{Wang2}  {\rm
Let $R$ be the two dimensional local Cohen-Macaulay ring
\[
k[\![X,Y,Z,U,V]\!]/(Z^2,ZU,ZV,UV,YZ-U^3,XZ-V^3),
\]
with $k$ a field and $X,Y,Z,U,V$ indeterminates. Let $I$ be the
maximal ideal ${\mathfrak m}$ of $R$. One has that the associated
graded ring ${\rm gr}_{{\mathfrak m}}(R)$ has depth zero and
\[
P_I(t)= \frac{1 + 3t + 3t^3 - t^4}{(1-t)^2}.
\]
In particular, one has $e_2 = e_1 - e_0 + 1$, that is, $e_2$ is
minimal according to Itoh's bound. However, $e_2$ is not minimal
with respect to the bound given in Theorem~\ref{thm1}.}
\end{Example}

Thus, the associated graded ring of the maximal ideal of the ring
$R$ can have depth zero even if $e_2 = e_1 - e_0 + 1$. More
generally, a condition such as $\lambda(R/I)=e_0-e_1+e_2$ is not
sufficient to guarantee that ${\rm gr}_I(R)$ is Cohen-Macaulay
even for an integrally closed ideal $I$. Motivated by this
failure, we observe that the right setting is the one of normal
ideals.

The following result is essentially contained in \cite{I}, we
present it for completeness with a simpler proof. As a piece of
notation, we denote by
$\overline{e}_0, \overline{e}_1, \ldots, \overline{e}_d$ the
Hilbert coefficients  with respect to the filtration ${\mathcal F} =
\{ \overline{I^n}
\}_{n\geq 0}$ given by the integral closure of the powers of $I$.

\begin{Theorem} \label{thm2}
Let $(R, {\mathfrak m})$ be a local Cohen-Macaulay ring of dimension
$d \geq 1.$ Let $I$
be an ideal generated by a system of parameters. If
$\lambda(R/\overline{I}) =
\overline{e}_0-\overline{e}_1+\overline{e}_2$ then
$\overline{I^{n+2}}=I^n\overline{I^2}$ for all $n \geq 0$.
\end{Theorem}
\demo Let $I=(x_1, \ldots, x_d)$. If $d=1, $ the result follows
from \cite[4.6]{HM}. Let assume then $d \geq 2.$  Write
$(\underline{x})=(x_1, \ldots, x_{d-2})$ and let $'$ denote images
modulo $(\underline{x})$. Call ${\mathfrak F} = \{
\overline{(I')^n} \}$. Then $\overline{e}_i = \overline{e}_i(I')$
for $i=0,1,2$ by \cite[8$($1$)$]{I}. Also by \cite[1$($2$)$]{I} we
have that $\lambda(R/\overline{I}) = \lambda(R'/\overline{I'})$.
This yields $\lambda(R'/\overline{I'}) =
\overline{e}_0(I')-\overline{e}_1(I')+\overline{e}_2(I')$. As
${\rm depth}\ {\rm gr}({\mathfrak F}) \geq 1 ={\rm dim}\,R'-1, $
by \cite[1.11$($4$)$]{GR} the degree of the $h$-polynomial
$s({\mathfrak F}) \leq 2$. Hence $\overline{e}_3(I')=0$. Now by
\cite[1.9]{GR} we obtain that $\overline{(I')^{n+1}} =
(x_{d-1},x_d)\overline{(I')^n}$ for $n\geq 2$, and, in particular,
$\overline{(I')^{n+2}} = (x_{d-1},x_d)^n\overline{(I')^2}$ for
$n\geq 0$. Finally, using \cite[17]{I}, we have that
$\overline{I^{n+2}}=(x_1, \ldots, x_d)^n\overline{I^2}$ for $n\geq
0$. \QED

\medskip

\begin{Theorem}\label{valla}
Let $(R, {\mathfrak m})$ be a local Cohen-Macaulay ring of
dimension $d \geq 1$. Let $I$ be an ${\mathfrak m}$-primary normal
ideal. The following conditions are equivalent:
\begin{itemize}
\item[$(${\it a}$)$]
$\lambda(R/I)=e_0-e_1+e_2$;

\item[$(${\it b}$)$]
$I^3=JI^2$ for some minimal reduction $J$ of $I$;

\item[$(${\it c}$)$]
$e_2 = \lambda(I^2/JI)$ for some minimal reduction $J$ of $I$.
\end{itemize}
Moreover, if any of the previous equivalent conditions holds then
${\rm gr}_I(R)$ is Cohen-Macaulay and
\[
P_I(t) = \frac{\lambda(R/I) + (e_0 - \lambda(R/I) -
\lambda(I^2/JI))t + \lambda(I^2/JI)t^2}{(1-t)^d}.
\]
\end{Theorem}
\demo Assume that condition $({\it a})$ holds and let $J=(x_1,
\ldots, x_d)$ be a minimal reduction of $I$. Observe that
${\mathcal{F}} = \{ \overline{J^n} \} = \{ I^n \}$, as $J^n$ is
still a reduction of $I^n$ $($not minimal$)$ hence
$\overline{e}_i(J) = e_i$. By Theorem~\ref{thm2}, applied to $J$,
we have $\overline{J^{n+2}} =J^n \overline{J^2}$ for $ n \geq 0$,
hence $I^{n+2}=J^nI^2$ for $n \geq 0$. This yields $I^3=JI^2$. Now
$({\it b})$ implies $({\it c})$ and $({\it c})$ forces $({\it a})$
by Theorem~\ref{thm1} and the Cohen-Macaulayness of ${\rm
gr}_I(R)$ and the Hilbert series follow as well from the same theorem. \QED

\section{Results on the higher Hilbert coefficients}

\noindent Itoh showed in \cite[3$($1$)$]{I} that $e_3 \geq 0$ for
an ${\mathfrak m}$-primary normal ideal. Earlier, Narita had given
in \cite {Nar} an example of a three dimensional Cohen-Macaulay
local ring and an ideal $I$ with $e_3 < 0$. The ring in Narita's
example contains nilpotents, thus Marley subsequently gave in
\cite[4.2]{M} an example of an ideal in a polynomial ring in three
variables with $e_3 <0$. Finally, the previously mentioned example
by Wang provides an example of a Cohen-Macaulay local ring $R$ in
which the maximal ideal has $e_3 <0$.

Let $n(I)$ denote the so-called postulation number of $I$, that is
the smallest integer $n$ such that $\lambda(R/I^n)$ is a
polynomial.

The following result improves the already known result of Itoh
under a weaker assumption.

\begin{Theorem}\label{thm3}
Let $(R, {\mathfrak m})$ be a Cohen-Macaulay local ring of
dimension three and with infinite residue field. Let $I$ be an
${\mathfrak m}$-primary ideal of $R$ such that $I^q$ is integrally
closed for some $q \geq n(I)$. Then $e_3 \geq 0$.
\end{Theorem}
\demo For $n \gg 0$ the Hilbert-Samuel function of $I$ can be
written as
\begin{equation}\label{In}
\lambda(R/I^n) = e_0 {n+2 \choose 3} - e_1 {n+1 \choose 2} + e_2
{n \choose 1} - e_3.
\end{equation}
Let $q \geq n(I)$ be an integer for which $\overline{I^q}=I^q$.
Consider the Hilbert-Samuel function of $I^q$. For $n\gg 0 $ one
has that
\begin{equation}\label{Inq}
\lambda(R/(I^q)^n) = \varepsilon_0 {n+2 \choose 3} - \varepsilon_1
{n+1 \choose 2} + \varepsilon_2 {n \choose 1} - \varepsilon_3.
\end{equation}
As $\lambda(R/(I^q)^n)=\lambda(R/I^{nq})$, an easy comparison
between $($\ref{In}$)$, with $nq$ in place of $n$, and
$($\ref{Inq}$)$ yields
\[
e_0 {nq+2 \choose 3} - e_1 {nq+1 \choose 2} + e_2 {nq \choose 1} -
e_3= \varepsilon_0 {n+2 \choose 3} - \varepsilon_1 {n+1 \choose 2}
+ \varepsilon_2 {n \choose 1} - \varepsilon_3
\]
or, equivalently,
\[
\frac{1}{6}e_0(n^3q^3+3n^2q^2+2nq)-\frac{1}{2}e_1(n^2q^2+nq)+e_2nq-e_3
= \qquad \qquad \qquad
\]
\[
\qquad \qquad \qquad \qquad
\frac{1}{6}\varepsilon_0(n^3+3n^2+2n)-\frac{1}{2}\varepsilon_1(n^2q+n)
+\varepsilon_2n-\varepsilon_3
\]
Hence one concludes that
\[
\varepsilon_0 = e_0 q^3 \qquad \varepsilon_1 = e_0 q^2(q-1) + e_1
q^2 \qquad \varepsilon_2 = e_0 {q \choose 3} + e_1 {q \choose 2} +
e_2 {q \choose 1} \qquad \varepsilon_3 = e_3.
\]
By \cite[12]{Itoh}, the Hilbert  coefficients of the ideal
$I^q$ satisfy the inequality
\[
\varepsilon_2 - \varepsilon_1 + \varepsilon_0 - \lambda(R/I^q)
\geq 0,
\]
as $q$ was chosen so that $\overline{I^q}=I^q$. After substituting
the $\varepsilon_i$'s with the corresponding expressions in terms
of the $e_i$'s we conclude that
\begin{eqnarray*}
\varepsilon_2 - \varepsilon_1 + \varepsilon_0 - \lambda(R/I^q) & =
& \left[e_0 {q \choose 3} + e_1 {q \choose 2} + e_2 {q \choose
1}\right] - \left[e_0 q^2(q-1) + e_1 q^2\right] + e_0 q^3 -
\lambda(R/I^q)
\\ & = & e_0 {q+2 \choose 3} - e_1 {q+1 \choose 2} + e_2 {q \choose
1} - \lambda(R/I^q) \\ & = & e_3,
\end{eqnarray*}
and therefore $e_3 \geq 0$, as claimed. \QED

\medskip

\begin{Remark}\label{ese} {\rm We note that in dimension three for
all $n\ge 2,$
$e_2(I^n) $ is always
strictly positive. In particular, the reduction number of $I^n$ is at
least two.}
\end{Remark}
\medskip

An ideal $I$ is said to be asymptotically normal if there exists
an integer $N \geq 1$ such $I^n$ is integrally closed for all $n
\geq N$. An interesting family of examples of asymptotically
normal ideals that are not normal are described  in the next
remark.

\medskip

\begin{Remark}\label{esempio}{\rm If $I$ is an asymptotically normal
ideal such that $e_3=0$ then ${\rm gr}_I(R)$ is not necessarily
forced to be Cohen-Macaulay. For example, the ideal $I$ in
Example~\ref{e2=0} is such that $I^n$ is integrally closed for
every $n \geq 2$, $e_3=0$ but ${\rm gr}_I(R)$ has depth zero.

More generally, let $I$ be any ${\mathfrak m}$-primary Gorenstein
ideal of minimal multiplicity in $R=k[\![X,Y,Z]\!]$, where $k$ is
a field and $X,Y,Z$ are indeterminates. Since $I \colon {\mathfrak
m}={\mathfrak m}^2, $ by \cite[3.6$({\it a})$]{CHV} we have that
${\mathfrak m}^2=I \colon {\mathfrak m} \subseteq \overline{I}$,
which forces $\overline{I}={\mathfrak m}^2 \neq I$.

Furthermore by \cite{Her}, $I/I^2 $ has a Cohen-Macaulay
deformation which is generically a complete intersection. Thus by
\cite[4.7.11, 4.7.17$({\it a})$]{BH} it follows that
$\lambda(I/I^2)=\lambda (R/I) \, {\rm ht}\, (I)=15$. Now
$\lambda(R/I^2)=\lambda (R/\mathfrak m^4) $ yields $I^2= \mathfrak
m^4$. By \cite[3.6$({\it b})$]{CHV} we also have that ${\mathfrak
m} I  = {\mathfrak m}^3$, which implies, in addition, that $I^3 =
I^2 I = {\mathfrak m}^4 I = {\mathfrak m}^3 {\mathfrak m} I =
{\mathfrak m}^6$. Hence we conclude that $I^n={\mathfrak m}^{2n}$
for all $n \geq 2$. Thus $I$ is asymptotically normal and its
Hilbert series is given by
\[
P_I(t)=\frac{5 + 6t^2 - 4t^3 + t^4}{(1-t)^3}.
\]
In particular $e_3=0$. On the other hand, ${\rm gr}_I(R)$ has
depth zero because for any superficial element $a \in I$ one has
$I^2 \colon a={\mathfrak m}^2 \not = I$.}
\end{Remark}

In the above remark, $I^2$ is a normal ideal in $k[\![X,Y,Z]\!] $
with $e_3(I^2)=0$, ${\rm gr}_{I^2}(R)$ Cohen-Macaulay and the
reduction number is two. It is natural to ask the following
question.

\begin{Question}\label{question}
Let $I$ be a normal ${\mathfrak m}$-primary ideal of a local
Gorenstein ring $R$. Does $e_3=0$ imply ${\rm gr}_{I }(R)$
Cohen-Macaulay? Does $e_3=0$ imply that the reduction number of
$I$ is two?
\end{Question}

In \cite[3]{Itoh} Itoh gave a positive answer to
Question~\ref{question} when $I$ is the maximal ideal of a
Gorenstein ring. Notice that if $I$ is asymptotically normal, but
not normal, the answer is negative by Remark~\ref{esempio}.

In Corollary~\ref{thm4} below we show that, in a local
Cohen-Macaulay ring of dimension $3$, the normality of $I$ implies
the Cohen-Macaulayness of ${\rm gr}_{I^n}(R)$ for all large $n$
whenever $e_3=0$.

\begin{Corollary}\label{thm4}
Let $(R, {\mathfrak m})$ be a local Cohen-Macaulay ring of
dimension three and with infinite residue field. Let $I$ be an
${\mathfrak m}$-primary ideal of $R$ such that $I$ is
asymptotically normal. Then $e_3=0$ if and only if there exists
some $n$ such that the reduction number of $I^n$ is at most two.
Under these circumstances then ${\rm gr}_{I^n}(R)$ is
Cohen-Macaulay.
\end{Corollary}
\demo Let $q \geq n(I)$ be an integer large enough so that $I^q$
is a normal ideal. From the proof of Theorem~\ref{thm3} we have
that
$0=e_3=\varepsilon_2-\varepsilon_1+\varepsilon_0-\lambda(R/I^q)$,
where the $\varepsilon_i$'s are the Hilbert coefficients of
$I^q$. The statement now follows from Theorem~\ref{valla}. \QED

\begin{Remark}{\rm A different proof of the above result can be
obtained in the following way. Let $q$ be the integer such that
$I^q$ is normal. By \cite[3.1]{HH}, there exists an integer $N$
such that ${\rm depth}\,{\rm gr}_{I^N}(R) \geq 2$ and by
\cite[4.6]{HM} we get
\[
e_3(I)=e_3(I^N)= \sum_{n\ge3} {{n-1}\choose 2}
\lambda(I^{nN}/JI^{nN-N}),
\]
for any minimal reduction $J$ of $I^N$. Hence $e_3(I)=0$ if and
only if $I^{3N}=JI^{2N}$ and the result follows.

The latter proof suggests the following result in dimension four.
By \cite[4.5, 4.1]{HM}, it is easy to obtain
\[
e_4(I)=e_4(I^N) \leq \sum_{n\ge4} {{n-1}\choose 3}
\lambda(I^{nN}/JI^{nN-N}),
\]
for any minimal reduction $J$ of $I^N$. }

\end{Remark}

\begin{Remark}{\rm
As we already remarked, a connection between the normality of $I$
and the depth of ${\rm gr}_{I^n}(R)$ has been observed in
\cite{HH}. Indeed Huckaba and Huneke show that if $I$ is normal
then ${\rm gr}_{I^n}(R)$ has depth at least $2$ for $n \gg 0$.
This result provides a two dimensional version of the
Grauert-Riemenschneider vanishing theorem. More precisely, this is
a generalization $($in dimension two$)$ of the following
formulation of Grauert-Riemenschneider due to Sancho de Salas: If
$R$ is a reduced Cohen-Macaulay local ring, essentially of finite
type over an algebraically closed field of characteristic zero,
and $I$ is an ideal of $R$ such that ${\rm Proj}(\mathcal R)$ is
regular, then ${\rm gr}_{I^n}(R)$ is Cohen-Macaulay for some $n
\gg 0$. While in dimension two the regularity of ${\rm
Proj}(\mathcal R)$ is not necessary $($as shown in \cite{HH}$)$,
in dimension three the Grauert-Riemenschneider theorem fails if
the assumption on ${\rm Proj}(\mathcal R)$ being regular is
dropped \cite{Cut}. In \cite{HH} Huckaba and Huneke give another
example of this failure.}
\end{Remark}

\begin{Example}
{\rm The same ideal $I$ considered in Example~\ref{ex-HuHu} also
shows that Corollary~\ref{thm4} is sharp, that is the condition on
$e_3$ cannot be relaxed. In fact, we checked that the ideal $I$ is
such that
\[
P_I(t)=\frac{31 + 43t +t^2 + t^3}{(1-t)^3}.
\]
Thus one has $e_3=1$. On the other hand, Huckaba and Huneke show
-- in \cite[3.11]{HH} -- that $I$ is a height $3$ normal $R$-ideal
for which ${\rm gr}_{I^n}(R)$ is not Cohen-Macaulay for any $n
\geq 1$.

In addition, one also has that $e_2=4$ while $\lambda(I^2/JI)=2$,
for any minimal reduction $J$ of $I$. Hence the bound in
Theorem~\ref{thm1} is strict in this setting.  }
\end{Example}

\end{document}